\documentclass[12pt]{article}
\usepackage[final]{epsfig}
\usepackage{graphics}
\usepackage{amsmath}
\usepackage{amsfonts}
\usepackage{latexsym}
\usepackage{amssymb}
\usepackage{graphicx}
\usepackage{url}
\usepackage{epstopdf}
\usepackage{hyperref}
\usepackage{color}
\usepackage{marginnote}
\usepackage{a4wide}

\newtheorem{lemma}{Lemma}[section]

\newtheorem{remark}[lemma]{Remark}

\newcommand{\g}{{\gamma}}

\newcommand{\proofend}{$\Box$\bigskip}

\newcommand{\R}{{\mathbb R}}
\newcommand{\Z}{{\mathbb Z}}

\def\proof{\paragraph{Proof.}}

\begin{document}

\title{Monotone twist maps and Dowker-type theorems}

\author{Peter Albers\footnote{
Mathematisches Institut,
Universit\"at Heidelberg,
69120 Heidelberg,
Germany;
peter.albers@uni-heidelberg.de}
 \and 
 Serge Tabachnikov\footnote{
Department of Mathematics,
Pennsylvania State University,
University Park, PA 16802,
USA;
tabachni@math.psu.edu}
} 

\date{\today}

\maketitle

\begin{abstract}
Given a planar oval, consider the maximal area of inscribed $n$-gons resp.~the minimal area of circumscribed $n$-gons. One obtains two sequences indexed by $n$, and one of Dowker's theorems states that the first sequence is concave and the second is convex. In total, there are four such classic results, concerning areas resp.~perimeters of inscribed resp.~circumscribed polygons, due to Dowker, Moln\'ar, and Eggleston. We show that these four results are all incarnations of the convexity property of Mather's $\beta$-function (the minimal average action function) of the respective billiard-type systems. We then derive new geometric inequalities of similar type for various other billiard system. Some of these billiards have been thoroughly studied, and some are novel. Moreover, we derive new inequalities (even for conventional billiards) for higher rotation numbers.
\end{abstract}


{\it \hfill To Misha Bialy with admiration on the occasion of his 60th birthday, belatedly.}

\section{Introduction} \label{sect:intro}

The classic Dowker theorem \cite{Do} concerns extremal polygons inscribed and circumscribed about an oval\footnote{Not to be confused with another classic Dowker theorem \cite{Do-other}.}. Here is its formulation.

\begin{figure}[ht]
\centering
\includegraphics[width=.7\textwidth]{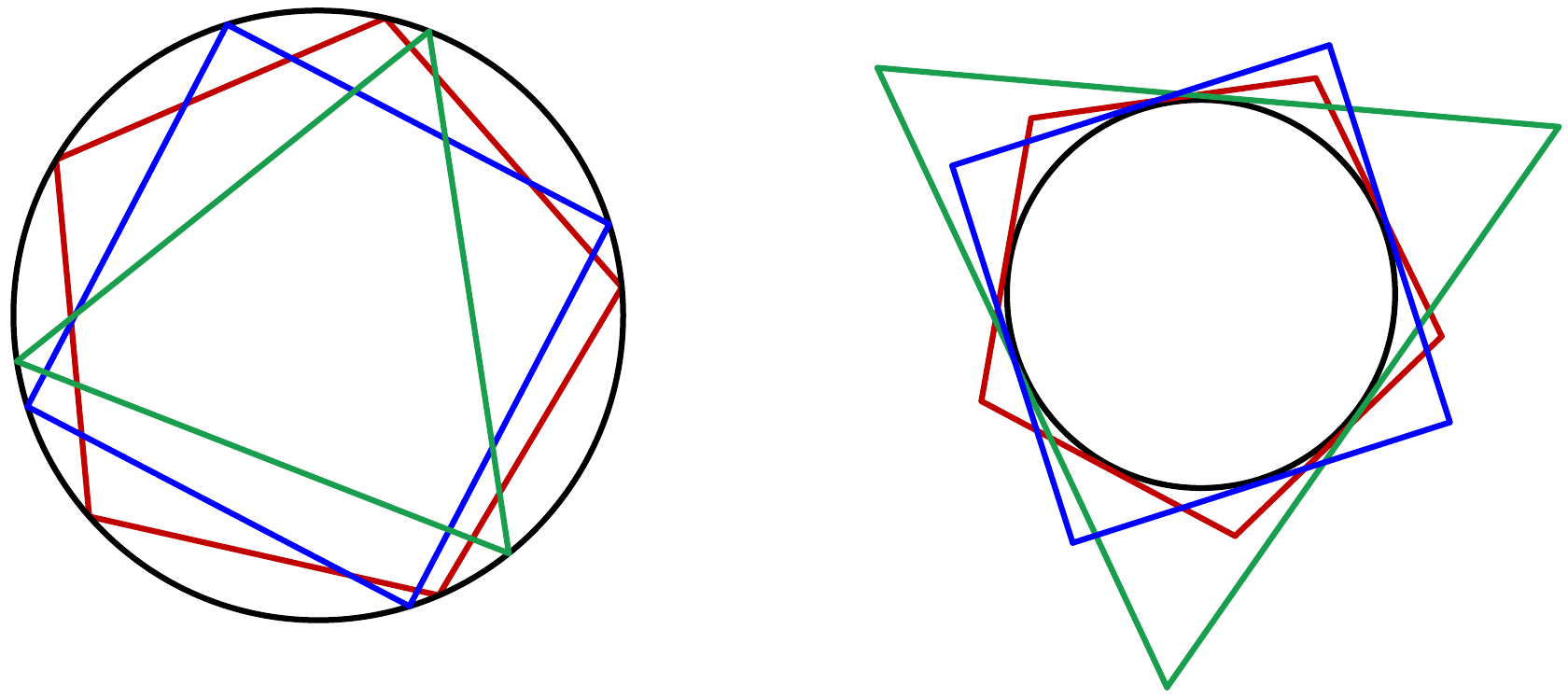}
\caption{Left: $P_{3}+P_{5}\le 2P_4$; right: $Q_{3}+Q_{5}\ge 2Q_4$.}	
\label{polys}
\end{figure}

Let $\g$ be a smooth strictly convex closed plane curve (an oval). Denote by $P_n$ the maximal area of $n$-gons inscribed in $\g$ and by $Q_n$ the minimal area of $n$-gons circumscribed about $\g$. Assume that $n\ge 4$. Then
\begin{equation} \label{eqn:Dow}
P_{n-1}+P_{n+1}\le 2P_n \quad {\rm and} \quad Q_{n-1}+Q_{n+1}\ge 2Q_n,
\end{equation}
see Figure \ref{polys}.
An analog of this result in the spherical geometry is due to L. Fejes T\'oth \cite{FT2}.

A similar result holds for perimeters. Let $R_n$ be the maximal perimeter of the $n$-gons inscribed in $\g$ and  $S_n$ be the minimal perimeter of the $n$-gons circumscribed about $\g$. A theorem, due to Moln\'ar \cite{Mo} and  Eggleston \cite{Eg}, states that
\begin{equation} \label{eqn:Mol}
R_{n-1}+R_{n+1}\le 2R_n\quad {\rm and} \quad S_{n-1}+S_{n+1}\ge 2S_n.
\end{equation}
See \cite{FT1,Gr,MSW,PA} for surveys and ramifications of these results.

In this article we show that these four results are particular cases of the convexity of the minimal average action function (Mather's $\beta$-function) of monotone twist maps, a result from the Aubry-Mather theory.  The maps in question are various kinds of billiards:  conventional billiards in $\g$ for the perimeter of inscribed polygons,  outer billiards about $\g$ for the area of circumscribed polygons,  symplectic billiards in $\g$ for the area of inscribed polygons, and  outer length billiards about $\g$ for the perimeter of circumscribed polygons.

The first two billiard systems have been thoroughly studied  for a long time; we refer to \cite{Ta} and \cite{DT} for surveys. Symplectic billiards were introduced only recently \cite{AT}. As to outer length billiards (the ``fourth billiards"), to the best of our knowledge, they were not studied before. We define them here and we plan to provide more details in the upcoming article \cite{AT2}.

We also apply the convexity of the minimal average action function to other billiard-like systems: wire billiards, wire symplectic billiards,  magnetic billiards, and outer magnetic billiards. 

Wire billiards were introduced and studied in \cite{BMT}; this is a dynamical system on the set of chords of a closed curve (satisfying certain conditions) in $\R^n$. The resulting Dowker-type inequality concerns the maximal perimeters of inscribed $n$-gons, that is, $n$-gons whose vertices lie on the curve.

Similarly, wire symplectic billiards is a dynamical system on the set of chords of a closed curve (also satisfying certain conditions) in linear symplectic space $\R^{2n}$. The resulting inequality concerns the maximal symplectic areas of inscribed $n$-gons. Wire symplectic billiards are introduced here for the first time. 

Magnetic billiards describe the motion of a charge in a magnetic field subject to the elastic reflections off the boundary of a plane domain. We consider the case of a weak constant magnetic field when the trajectory of a charge comprises circular arcs of a fixed radius that is greater than the greatest radius of curvature of the boundary of the domain, the oval $\g$. The resulting Dowker-style geometric inequality combines the perimeter of a trajectory with the area bounded by it.  

Outer magnetic billiards are similar to outer billiards, but instead of tangent lines one considers tangent arcs of a sufficiently great fixed radius. The resulting inequality concerns the minimal areas of circumscribed curvilinear $n$-gons. Outer magnetic billiards were introduced and studied in \cite{BM}. 

Let us emphasize that the classic Dowker-type theorems concern simple polygons, whereas our inequalities include star-shaped polygons as well: the former have rotation numbers $1/q$, and the latter have the more general rotation numbers $p/q$. The corresponding inequalities are new even in the area/perimeter and inscribed/circumscribed cases.

We hope that the topic of this article is of interest to two research communities, the Hamiltonian dynamical and the convex geometrical ones. They use different methods to study closely related problems, and an interaction of these two communities would be beneficial.  We make some comments   on  applications of the theory of interpolating Hamiltonians in convex geometry at the end of the article. 

\section{Monotone twist maps and minimal action} \label{sect:ineq}

Let us recall  basic facts about monotone twist maps; we refer to \cite{CMSS,Go,KH, Me,Si}. 

We consider a cylinder $S^1\times (a,b)$ where $-\infty\leq a<b\leq+\infty$ and an area preserving diffeomorphism $f:S^1\times (a,b)\to S^1\times (a,b)$, isotopic to the identiy. Let $F=(F_1,F_2):\R\times (a,b)\to\R\times (a,b)$ be a lift to the universal cover.
We denote by $(x,y)\in\R\times (a,b)$ the coordinates on the strip $\R\times (a,b)$. In these coordinates the area form is $dy\wedge dx$. The monotone twist condition is 
\begin{equation*}\label{eqn:twist_condition_1}
\frac{\partial F_1(x,y)}{\partial y} > 0.
\end{equation*}
If $a$, resp.~$b$, is finite, we assume that  $F$ extends to the boundary of the strip as a rigid shift, that is $F(x,a)=(x+\omega_-,a)$, resp.~$F(x,b)=(x+\omega_+,b)$. Otherwise, we set $\omega_-:=-\infty$, resp.~$\omega_+:=+\infty$. We call $f$ a monotone twist map. The interval $(\omega_-,\omega_+)\subset\R$ is called the twist interval of the map $f$. The twist condition and the twist interval do not depend on the choice of the lift $F$. 

In the following sections $S^1\times(a,b)$ will appear as (diffeomorphic copy of) the phase space of various billiard systems. The foliation of the phase space corresponding to  $\{\mathrm{pt}\}\times (a,b)$ will be referred to as vertical foliation. Moreover, for now, $S^1=\R/\Z$ for simplicity, but in the billiard situations $S^1$ has varying length, e.g.,~for conventional billiards the length is the perimeter of the billiard table. 

Monotone twist maps can be defined via generating functions. A function 
$$
H: \{(x,x')\in\R\times \R\mid \omega_-<x'-x<\omega_+\} \longrightarrow \R
$$ 
is a generating function for $f$ if the following holds:
\begin{equation}\label{eqn:generating_function_1}
F(x,y)=(x',y') \quad\text{if and only if}\quad \frac{\partial H(x,x')}{\partial x}=-y,\ \frac{\partial H(x,x')}{\partial x'}=y'.
\end{equation}
The variables are related via the diffeomorphism $(x,y)\mapsto (x,x')=(x,F_1(x,y))$. The function $H$ is periodic in the diagonal direction: $H(x+k,x'+k)=H(x,x')$ for $k\in\Z$.
The twist condition becomes the inequality
\begin{equation}\label{eqn:twist_condition_2}
\frac{\partial^2 H}{\partial x\partial x'}(x,x')<0.
\end{equation}
In the coordinates $(x,x')$ condition \eqref{eqn:generating_function_1} becomes
\begin{equation}\label{eqn:generating_function_2}
F(x,x')=(x',x'') \quad\text{if and only if}\quad \frac{\partial }{\partial x'} \Big(H(x,x')+H(x',x'')\Big)=0.
\end{equation}
As a consequence, the differential 2-form 
$$
\frac{\partial^2 H}{\partial x\partial x'}(x,x')\ dx\wedge dx'
$$
is invariant under the map.

In terms of the coordinates $(x,x')$, a generating function $H(x,x')$ is not uniquely defined, for instance it can be changed to $H(x,x')+h(x')-h(x)$ by any function $h(x)$ such that $h(x+1)-h(x)$ is constant. This changes the coordinates $y$ and $y'$, but does not change the periodicity in the diagonal directions of the function $H(x,x')$, the twist condition \eqref{eqn:twist_condition_2}, or the variational characterization \eqref{eqn:generating_function_2}. In terms of $(x,y)$-coordinates the map $F$ is conjugated by the area-preserving diffeomorphism $(x,y)\mapsto(x,y+h'(x))$. This explains for instance the different conventions for generating functions in the literature, e.g.,~for conventional billiards. 

Birkhoff periodic orbits for $f$ of type $(p,q)\in\Z\times\Z_+$ are  bi-infinite sequences $(x_n,y_n)_{n\in \Z}$ in $\R\times (a,b)$ such that for all $n\in \Z$ we have
\begin{enumerate}
\item[1)] $x_{n+1} > x_n$;
\item[2)] $(x_{n+q},y_{n+q})=(x_n+p,y_n)$;
\item[3)] $F(x_n,y_n)=(x_{n+1},y_{n+1})$.
\end{enumerate}
The Birkhoff theorem asserts that for every rational number $\tfrac{p}{q} \in (\omega_-,\omega_+)$ in lowest terms, the map $f$ possesses 
at least two Birkhoff periodic orbits of type $(p,q)$. One of these orbits (the ``easier one") corresponds to the minimum of the function
\begin{equation} \label{eqn:funct}
H(x_0,x_1)+H(x_1,x_2)+\ldots + H(x_{q-1},x_q)
\end{equation}
on the space of bi-infinite sequences of real numbers $X=(x_n)_{n\in\Z}$ satisfying the monotonicity condition $x_{n+1} \ge x_n$ and the periodicity condition  $x_{n+q}=x_n+p$. Setting $y_n:=\frac{\partial H}{\partial x'}(x_{n-1},x_n)$, it turns out that, due to the twist condition \eqref{eqn:twist_condition_2}, the sequence $(x_n,y_n)$ actually satisfies the above conditions 1) -- 3), in particular the stronger monotonicity condition 1).

Let us denote the minimal value of the function \eqref{eqn:funct} on this space of bi-infinite sequences by $T_{p,q}$. Then the minimal action of a $(p,q)$-periodic orbit is defined as follows:
\begin{equation}\label{eqn:def_beta_fctn}
\beta\left(\frac{p}{q}\right):= \frac{1}{q}\ T_{p,q}.
\end{equation}
This is the celebrated Mather $\beta$-function. The amazing fact is that Mather's $\beta$-function is well-defined, i.e., does not depend on the representation of the rational number $\frac{p}{q}$, for example, $\beta(\frac26)=\frac16T_{2,6}=\frac13T_{1,3}=\beta(\frac13)$.

We note that the above mentioned change of the generating function $\tilde H(x,x')=H(x,x')+h(x')-h(x)$ with $h(x+1)-h(x)=c$ leads to the change of the $\beta$-function:
$$
\tilde \beta\left(\frac{p}{q}\right) = \beta\left(\frac{p}{q}\right) + c\, \frac{p}{q}.
$$

Mather's $\beta$-function is a strictly convex continuous function of the rotation number: if $\frac{p_1}{q_1} \neq \frac{p_2}{q_2}$, then 
\begin{equation} \label{eqn:convex}
t \beta\left(\frac{p_1}{q_1}\right) + (1-t) \beta\left(\frac{p_2}{q_2}\right) > \beta\left(t\ \frac{p_1}{q_1}+ (1-t)\ \frac{p_2}{q_2}\right)
\end{equation}
for all $t\in(0,1)$, see \cite{MF,Si}. Although the minimal average action extends as a strictly convex function to irrational rotations numbers too, we only consider rational ones, therefore, in what follows, $t$ is also a rational number. 
The following lemma deduces the general Dowker-style inequality from the convexity of the minimal action.

\begin{lemma} \label{lm:ineq}
For all relatively prime $(p,q)\in\Z\times\Z_+$ with $\tfrac{p}{q} \in (\omega_-,\omega_+)$ and $q\neq1$, the inequality
$$
T_{p,q-1}+T_{p,q+1} > 2T_{p,q}
$$
holds.
\end{lemma} 

\proof
Consider the inequality \eqref{eqn:convex} with the choices 
$$
p_1=p_2=p,\ q_1=q-1,\ q_2=q+1,\  t=\frac{q-1}{2q},
$$
i.e.,
\begin{equation}\nonumber
\tfrac{q-1}{2q}\beta(\tfrac{p}{q-1})+\big(1-\tfrac{q-1}{2q}\big)\beta(\tfrac{p}{q+1})>\beta\Big(\tfrac{q-1}{2q}\tfrac{p}{q-1}+\big(1-\tfrac{q-1}{2q}\big)\tfrac{p}{q+1}\Big),
\end{equation}
and simplify to
\begin{equation}\nonumber
\begin{aligned}
\tfrac{q-1}{2q}\beta(\tfrac{p}{q-1})+\tfrac{q+1}{2q}\beta(\tfrac{p}{q+1})&>\beta\Big(\tfrac{q-1}{2q}\tfrac{p}{q-1}+\tfrac{q+1}{2q}\tfrac{p}{q+1}\Big)\\
&=\beta\Big(\tfrac{p}{2q}+\tfrac{p}{2q}\Big)
=\beta\big(\tfrac{p}{q}\big).
\end{aligned}
\end{equation}
Then \eqref{eqn:def_beta_fctn} gives
\begin{equation}\nonumber
\tfrac{q-1}{2q}\tfrac{1}{q-1}T_{p,q-1}+\tfrac{q+1}{2q}\tfrac{1}{q+1}T_{p,q+1}=\tfrac{q-1}{2q}\beta(\tfrac{p}{q-1})+\tfrac{q+1}{2q}\beta(\tfrac{p}{q+1})>\beta\big(\tfrac{p}{q}\big)=\tfrac{1}{q}T_{p,q},
\end{equation}
which simplifies to
\begin{equation}\nonumber
\tfrac{1}{2q}T_{p,q-1}+\tfrac{1}{2q}T_{p,q+1}>\tfrac{1}{q}T_{p,q},
\end{equation}
i.e.,
\begin{equation}\nonumber
T_{p,q-1}+T_{p,q+1}>2T_{p,q},
\end{equation}
as claimed.
\proofend

%
%

Next we describe a small extension of the above discussion, well-known to experts. It is sometimes convenient to consider as phase space a set of the form
$$
\{(\bar x,\bar y)\mid \bar x\in \R,\; 0\leq \bar y\leq o(\bar x)\}
$$
for some function $o:\R\to(0,\infty)$ together with an area-preserving (with respect to $d\bar y\wedge d\bar x)$ self-map $\bar F$ which is a monotone twist map, i.e.,~$\frac{\partial \bar F_1}{\partial \bar y}>0$. This set-up can be transformed to the above standard setting. First, we observe that the map
\begin{equation}\label{eqn:area_preserving_strip_to_wavy}
\begin{aligned}
\R\times[0,1]&\to\{(\bar x,\bar y)\mid \bar x\in \R,\; 0\leq \bar y\leq o(\bar x)\}\\
(x,y)&\mapsto \big(O^{-1}(x), o(O^{-1}(x))y\big)
\end{aligned}
\end{equation}
is an area-preserving diffeomorphism if $O(\bar x)$ is an anti-derivative of $o(\bar x)$. Note that, since $O'(\bar x)=o(\bar x)>0$, the function $O$ is strictly monotone and, in particular, invertible.

For simplicity, set $\varphi(x):=O^{-1}(x)$. Then the above maps reads
\begin{equation}\nonumber
\begin{aligned}
\bar x&=\varphi(x)\\
\bar y &= o(\varphi(x))y=o(\bar x)y\;.
\end{aligned}
\end{equation}
That is, $\R\times\{0\}$ is mapped to itself and $\R\times\{1\}$ to $\{\bar y=o(\bar x)\}$. The map \eqref{eqn:area_preserving_strip_to_wavy} is area-preserving since
\begin{equation}\nonumber
d\bar x\wedge d\bar y= \varphi'(x)o(\varphi(x))dx\wedge dy
\end{equation}
and
\begin{equation}\nonumber
\varphi'(x)o(\varphi(x))=\frac{d}{dx}O(\varphi(x))=1
\end{equation}
since $\varphi(x)=O^{-1}(x)$. 

We point out that the two vertical foliations (given by fixing $x$, resp. $\bar x$) are mapped to each other by the map \eqref{eqn:area_preserving_strip_to_wavy}. Moreover, if we have area-preserving maps $\bar F$ on $\{(\bar x,\bar y)\mid \bar x\in \R,\; 0\leq \bar y\leq o(\bar x)\}$ and $F$ on $\R\times[0,1]$ which are conjugate to each other by \eqref{eqn:area_preserving_strip_to_wavy}, then one satisfies the twist condition if and only if the other does. This uses again $o(x)>0$. Finally, the twist condition and the variational description of $\bar F$ in terms of a generating function continues to hold.

\section{The  classic Dowker-type theorems revisited} \label{sect:old}

In this section we consider four billiard-like systems: two inner, two outer, two with  length, and two with  area, as their generating functions. Lemma \ref{lm:ineq} will imply the four  Dowker-style theorems mentioned in the introduction.

\subsection{Conventional  billiards} \label{subsect:conv}

 Let $\g:S^1\to\R^2$ be a closed smooth strictly convex planar curve (an oval), oriented counterclockwise and parameterized by arc length. This is the boundary of the billiard table. We assume that it has unit length. Let $x$ be the respective coordinate on $\R$, the universal cover of $S^1$.
 
The phase space of the billiard ball map $F$ is the space of oriented chords of $\g$, the vertical foliation consists of the chords with a fixed initial point. The generating function is given by the formula\footnote{There are various conventions in the literature. Another common choice of generating function is $-|\g(x) \g(x')|$. As explained in the previous section, adding $x'-x$ does not change the twist condition, etc.}
$$
H(x,x')=x'-x-|\g(x) \g(x')|,
$$
where $|\g(x) \g(x')|$ denotes the chord length, i.e.,~the Euclidean length of the segment between $\g(x)$ and  $\g(x')$. 

One can calculate (see, e.g., \cite{Ta}) that
$$
\frac{\partial |\g(x) \g(x')|}{\partial x} = -\cos\alpha,\ 
\frac{\partial |\g(x) \g(x')|}{\partial x'} = \cos\alpha',
\ {\rm and}\ \ 
\frac{\partial^2 |\g(x) \g(x')|}{\partial x \partial x'} = \frac{\sin \alpha \sin \alpha'}{|\g(x) \g(x')|},
$$
where $\alpha,\alpha'\in(0,\pi)$ are the angles made by the chord $\g(x) \g(x')$ with the curve $\g$. Therefore, 
as coordinates we obtain
$$
y=1-\cos\alpha,\ y'=1-\cos\alpha'\in(0,2),
$$
and
$$
\frac{\partial^2 H(x,x')}{\partial x \partial x'} = -\frac{\sin \alpha \sin \alpha'}{|\g(x) \g(x')|}<0.
$$ 
It follows that the quantity $T_{p,q}$ from section \ref{sect:ineq} is $p$ minus the greatest perimeter of the  $q$-gons with the rotation number $p$ inscribed in $\g$.  Denoting this perimeter by $R_{p,q}$,  Lemma \ref{lm:ineq} implies that $R_{p,q-1}+R_{p,q+1} < 2R_{p,q}$ which, for $p=1$, reduces to the statement of the Moln\'ar-Eggleston theorem (\ref{eqn:Mol}). 

\begin{figure}[ht]
\centering
\includegraphics[width=.3\textwidth]{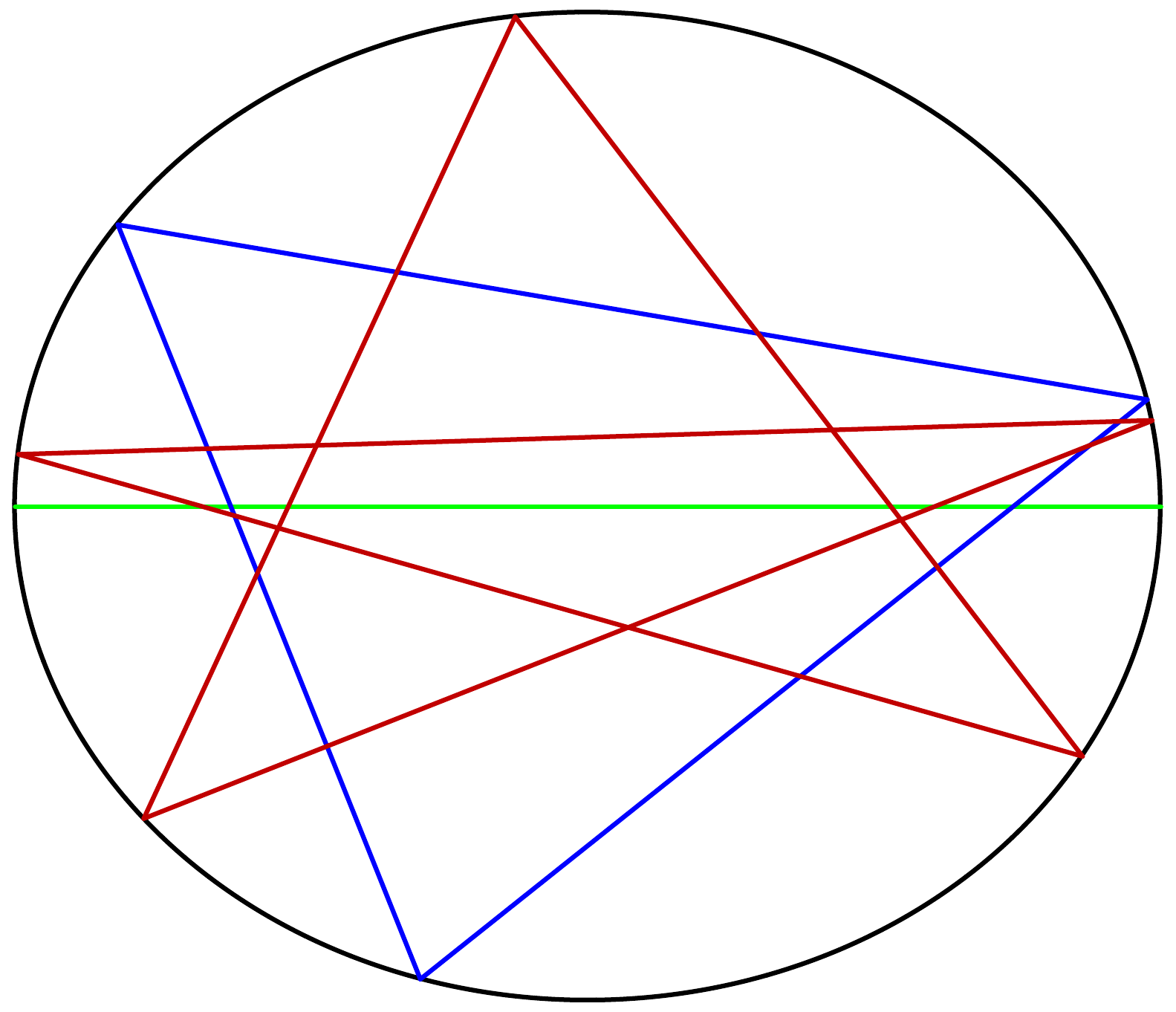}
\caption{The curve is an ellipse with the eccentricity of about $0.5$. One has  $\frac{R_{2,4} + R_{2,6}}{2 R_{2,5}}\approx 99\%$.}	
\label{456}
\end{figure}

As explained above, even for conventional billiards the inequalities $R_{p,q-1}+R_{p,q+1} < 2R_{p,q}$ are new for $p>1$. In Figure \ref{456} we illustrate the inequality $R_{2,4} + R_{2,6} < 2 R_{2,5}$. The extreme quadrilateral with $p=2$ is the diameter of the curve, traversed four times, i.e., $R_{2,4}=2R_{1,2}$, and the extreme hexagon with $p=2$ is the extreme triangle, traversed twice, i.e., $R_{2,6}=2R_{1,3}$. Therefore, the inequality can be rewritten as $R_{1,2} + R_{1,3} < R_{2,5}$, which has a different form from the Dowker-style inequalities. 

\subsection{Outer (area) billiards} \label{subsect:out}

\begin{figure}[ht]
\centering
\includegraphics[width=.39\textwidth]{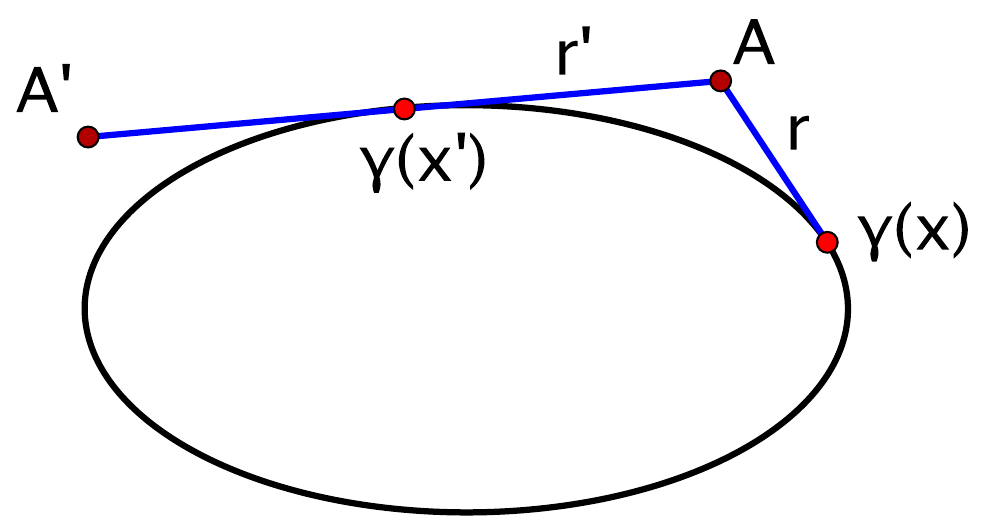}
\caption{The outer billiard map.}	
\label{outer}
\end{figure}

The outer billiard map $F$ about a smooth closed oriented strictly convex curve $\g$ is depicted in Figure \ref{outer}: one has $F(A)=A'$ if the orientation of the segment $AA'$ coincides with the orientation of $\g$, and $|Ax'|=|x'A'|$, see \cite{DT} for a survey. In this article we sometimes call this dynamical system outer area billiards to distinguish it from the upcoming outer length billiard in Section \ref{subsect:outlength}.

The map $F$ is an area preserving map of the exterior of $\g$ with respect to the standard area of the plane. It extends as the identity to $\g$. The phase space is foliated by the positive tangent rays to $\g$, and $F$ is a twist map.
Let $x\in\R$ be the (lifted to $\R$) angular coordinate on $\g$, that is, the direction of the oriented tangent line of $\g$. Moreover, if we write a point $A$ in the exterior of $\g$ as $\g(x)+r\g'(x)$, $r>0$, then the map 
$$\{\text{exterior of }\g\}\ni A\mapsto (x,r)\in S^1\times(0,\infty)$$
 is a symplectomorphism between the standard area form and $rdr\wedge dx=dy\wedge dx$ with $y=\tfrac{r^2}{2}$.

The generating function $H:\{(x,x')\mid 0< x'-x<\pi\}\to\R$ of the map $F$ is the area of the (oriented) curvilinear triangle obtained by first following the segments $\g(x)A$ and $A\g(x')$ and then the arc $\g(x')\g(x)$, see Figure \ref{outer}.
One has
$$
\frac{\partial H(x,x')}{\partial x} = -\frac{r^2}{2},\quad \frac{\partial H(x,x')}{\partial x'} = \frac{(r')^2}{2},
$$
and $$
\frac{\partial^2 H(x,x')}{\partial x \partial x'} =-r\frac{\partial r}{\partial x'}<0.
$$
It follows that the quantity $T_{p,q}$ is the minimal area of the circumscribed $q$-gon with the rotation number $p$ minus a constant ($p$ times the area bounded by $\g$). Denoting this circumscribed area by $Q_{p,q}$, 
Lemma \ref{lm:ineq} implies that $Q_{p,q-1}+Q_{p,q+1} > 2Q_{p,q}$ which, for $p=1$, is a statement of the Dowker theorem (\ref{eqn:Dow}). 

\subsection{Symplectic billiards} \label{subsect:sympl}

Symplectic billiards were introduced and studied in \cite{AT}, see also the recent papers \cite{BB,BBN}. 

Let $\g(x)$ be a positively oriented parameterized smooth closed strictly convex planar curve. For a point $\g(x)$, let $\g(x^*)$ be the other point on $\g$ where the tangent line is parallel to that in $\g(x)$. The phase space of symplectic billiard is then the set of the oriented chords $\g(x)\g(x')$ where $x<x'<x^*$ according to the orientation of $\g$. That is, the phase space is the set of pairs $(x,x')$ such that $\omega(\g'(x),\g'(x')) > 0$. Here, $\omega$ is the standard area form in the plane, the determinant made by two vectors. 

\begin{figure}[ht]
\centering
\includegraphics[width=.3\textwidth]{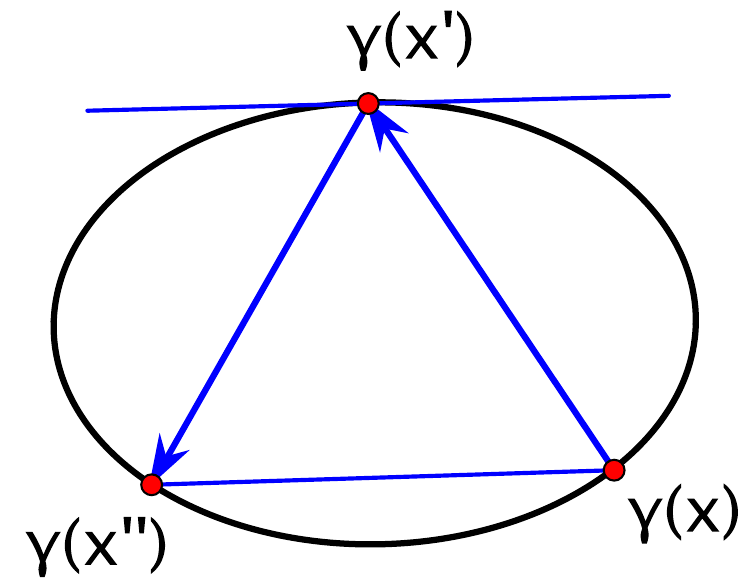}
\caption{The symplectic billiard map.}	
\label{symplectic}
\end{figure}

The vertical foliation consists of the chords with a fixed initial point. The symplectic billiard map $F$ sends a chord $\g(x)\g(x')$ to $\g(x')\g(x'')$ if the tangent line $T_{\g(x')} \g$ is parallel to the line $\g(x)\g(x'')$, see Figure \ref{symplectic}. Unlike the conventional billiards, this reflection law is not local. We note that if $\omega(\g'(x),\g'(x')) > 0$, then $\omega(\g'(x'),\g'(x'')) > 0$ as well (see \cite{AT}). 

We extend the map $F$ to the boundary of the phase space by continuity: $F(x,x):=(x,x)$ and $F(x,x^*):=(x^*,x)$.

The generating function $H:\{(x,x')\mid x<x'<x^*\}\to\R$ is the area bounded by the oriented bi-gon formed by following first the arc $\g(x)\g(x')$ and then the segment $\g(x')\g(x)$.
Note the similarity of this generating function with the one of the conventional  billiard: the length is replaced by the area. One has
$$
y=-\frac{\partial H(x,x')}{\partial x}=\frac12\omega(\g'(x),\g(x')-\g(x))
$$
and
$$
\frac{\partial^2 H(x,x')}{\partial x \partial x'} = -\frac{1}{2} \omega(\g'(x),\g'(x')) < 0.
$$ 
Thus, this is a situation where the phase space naturally is of the form 
$$
\{(\bar x,\bar y)\mid \bar x\in \R,\; 0\leq \bar y\leq o(\bar x)\}
$$
where $o(x)=\frac12\omega(\g'(x),\g(x^*)-\g(x))$.

Let $\g(x),\g(x'),\g(x'')$ be three consecutive vertices of an inscribed polygon of the maximal area. Then $T_{\g(x')} \g$ is parallel to the line $\g(x)\g(x'')$, hence $F(x,x')=(x',x'')$.  Then either $\omega(\g'(x),\g'(x')) > 0$, or $\omega(\g'(x),\g'(x')) < 0$, or $\omega(\g'(x),\g'(x')) = 0$. One has the same inequality for $\omega(\g'(x'),\g'(x''))$. If this is negative, the total area is negative as well, and changing the orientation makes it greater. Likewise if $\omega(\g'(x),\g'(x')) = 0$, then the area is not maximal either. Hence the maximal area polygon corresponds to the Birkhoff minimal $(p,q)$-periodic orbit of the symplectic billiard map. 

It follows that, up to an additive constant ($p$ times the area bounded by $\g$), the quantity $T_{p,q}$ is minus the greatest area of the  $q$-gon with the rotation number $p$, inscribed in $\g$.
Denoting this area by $P_{p,q}$,  Lemma \ref{lm:ineq} implies that $P_{p,q-1}+P_{p,q+1} < 2P_{p,q}$ which, for $p=1$, is a statement of the Dowker theorem (\ref{eqn:Dow}). 

\subsection{Outer length  billiards} \label{subsect:outlength}

As far as we know, the outer length billiard system, defined by extremizing the perimeter of a circumscribed polygons, has not been described in the literature yet. We provide necessary details here, and will return to a more detailed study of this dynamical system in \cite{AT2}. See \cite{DeT} for a study of the  polygons circumscribed about a convex curve and having the minimal perimeter.

\begin{figure}[ht]
\centering
\includegraphics[width=.4\textwidth]{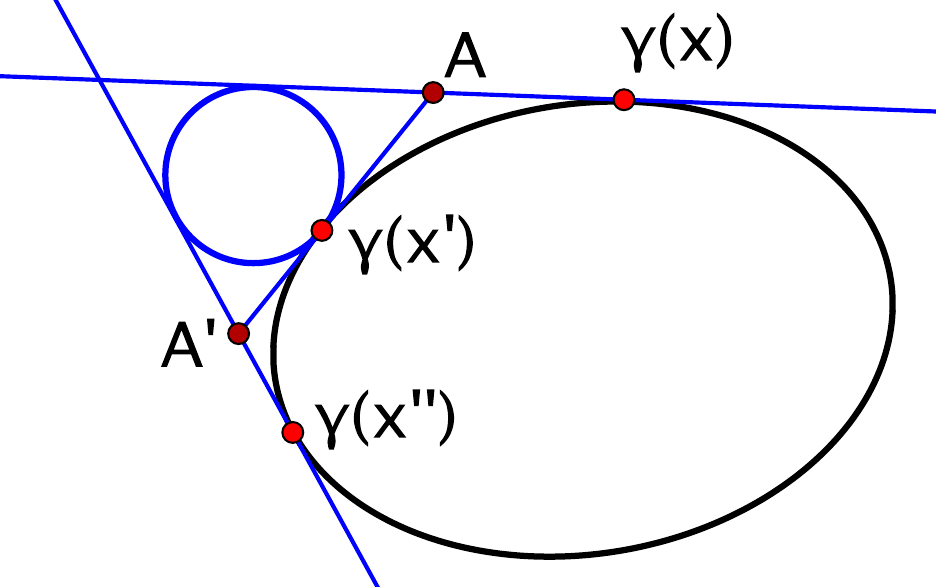}
\caption{Outer length billiard.}	
\label{fourth}
\end{figure}

The map $F$ acts on the exterior of an arc length parameterized  oval $\g(x)$ and is given by the following geometrical construction, see Figure \ref{fourth}. Let $A$ be a point outside of  $\g$, and let $A\g(x')$ and $A\g(x)$ be the positive and negative tangent segments to $\g$ (the sign given by the orientation of the oval). Consider the circle tangent to $\g$ at the point $\g(x')$, tangent to the line $A\g(x)$, and lying on the right of the ray $A\g(x')$. Then $A'=F(A)$ is defined as the intersection point of the line $A\g(x')$ and the common tangent line of the circle and $\g$ (tangent to $\g$ in $\g(x'')$ in Figure \ref{fourth}).

The map $F$ extends as identity to $\g$. The vertical foliation in phase space, i.e., the exterior of $\g$, is the same as for the outer (area) billiard: it consists of the positive tangent rays to $\g$. Similarly to outer billiards, $F$ is a twist map.  We recall from symplectic billiards  that the points $\g(x)$ and $\g(x^*)$ have parallel tangent lines.

\begin{lemma} \label{lm:genout}
The generating function $H:\{(x,x')\mid x<x'<x^*\}\to\R$ of the map $F$ is given by the formula 
$$
H(x,x')=|\g(x)A|+|A\g(x')|-x'+x.
$$ 
\end{lemma}

\proof
Fix points $\g(x)$ and $\g(x'')$ and consider
\begin{equation}\nonumber
\begin{aligned}
H(x,x')+H(x',x'')&=|\g(x)A|+|A\g(x')|+|\g(x')A'|+|A'\g(x'')|-x''+x'-x'+x\\
&=|\g(x)A|+|AA'|+|A\g(x'')|-x''+x.
\end{aligned}
\end{equation}
We claim that the specific point $\g(x')$ described above, see Figure \ref{fourth}, extremizes the length $|\g(x)A|+|AA'|+|A'\g(x'')|$. 

\begin{figure}[ht]
\centering
\includegraphics[width=.47\textwidth]{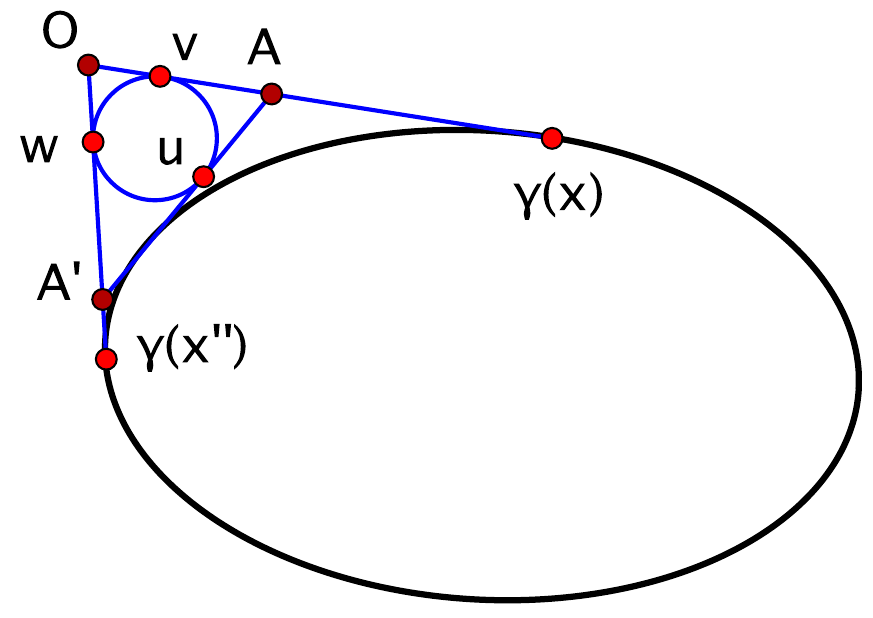}
\caption{Construction of the extremizer.}	
\label{fourpf}
\end{figure}

To prove this claim, consider Figure \ref{fourpf}. We use the fact that the two tangent segments to a circle through a common point have equal lengths, e.g.~$|Au|=|Av|$ in Figure \ref{fourpf}. 

Therefore, we can rewrite as follows
$$
|\g(x)A|+|AA'|+|A'\g(x'')|=|\g(x)v|+|\g(x'')w|=|\g(x)O|+|\g(x'')O|-2|Ov|.
$$
The left hand side is minimal when $|Ov|$ is maximal, and this happens when the circle is greatest possible, i.e., if $u$ lies on $\g$, that is, $u=\g(x')$ in Figure \ref{fourth}.

A similar argument applies if the point $O$ is on the other side of the line $AA'$, and if the tangent lines at points $\g(x)$ and $\g(x'')$ are parallel.
\proofend

The partial derivatives of the generating function are
$$
y=-\frac{\partial H(x,x')}{\partial x}= k(x)|A\g(x)| \cot \frac{\varphi}{2}\in(0,\infty),
$$
where $k$ is the curvature function of $\g$, and $\varphi$ is the angle between the tangent segments $A\g(x)$ and $A\g(x')$, and 
$$
\frac{\partial^2 H(x,x')}{\partial x\partial x'}= - \frac{k(x)k(x')(|A\g(x)|+|A\g(x')|)}{2\sin^2 \frac{\varphi}{2}}.
$$
It follows that, up to an additive constant ($p$ times the perimeter of $\g$), the quantity $T_{p,q}$ is  the minimal perimeter of the  $q$-gons with the rotation number $p$, circumscribed about $\g$.  Denoting this perimeter by $S_{p,q}$,  Lemma \ref{lm:ineq} implies that $S_{p,q-1}+S_{p,q+1} > 2S_{p,q}$ which, for $p=1$, is a statement of the Moln\'ar-Eggleston theorem (\ref{eqn:Mol}).

\begin{remark} \label{rmk:area}
{\rm
The area form that is invariant under this billiard map is, in terms of the generating function, 
$$
-\frac{\partial^2 H(x,x')}{\partial x \partial x'}\ dx\wedge dx'.
$$
This is a functional multiple of the standard area form $\omega$ in the exterior of the oval and,
at the point $A$ (in Figure \ref{fourth}), its value is, see \cite{AT2},
$$
\cot \left(\frac{\varphi}{2}  \right) \left( \frac{1}{|A\g(x)|}+\frac{1}{|A\g(x')|} \right) \omega.
$$
}
\end{remark}

\begin{remark} \label{rmk:Laz}
{\rm
The quantity $H=|\g(x)A|+|A\g(x')|-x'+x$ is known in the study of (the conventional) billiards as the Lazutkin parameter. Given an oval $\g$, consider the locus of points $A$ for which $H$ has a constant value. This locus is a curve $\Gamma$, and the billiard inside $\Gamma$ has the curve $\g$ as a caustic: a billiard trajectory tangent to $\g$ remains tangent to it after the reflection in $\Gamma$. This is known as the string construction of a billiard curve by its caustic (see, e.g., \cite{Ta}).

A similar relation exists between the level curves of the generating function of the symplectic billiard and the invariant curves of the outer billiard. Consider the set of chords that cut off a fixed area from an oval $\g$, that is,  a level curve of the generating function of the symplectic billiard in $\g$. The envelope of these chords is a curve $\Gamma$, and $\g$ is an invariant curve of the outer billiard about $\Gamma$. This is the area construction of an outer billiard curve by its invariant curve. 

The meaning of this relation between the level curves of a generating function of one billiard system and invariant curves of another one is not clear to us.
}
\end{remark}

\section{More examples} \label{sect:ex}

\subsection{Inner and outer billiards in $S^2$ and $H^2$} \label{subsect:curved}

 Inner (conventional) billiards are defined in the same way in the spherical and hyperbolic geometries as in the Euclidean plane: the boundary of the billiard table $\g$ has positive geodesic curvature and, in the case of $S^2$, this implies that $\g$ is contained in an open hemisphere. The billiard ball travels along geodesics and reflects in $\g$ subject to the law of equal angles.

Outer billiards are defined in $H^2$ similarly to the Euclidean case, but the case of $S^2$ is somewhat different.

Let $\g$ be a closed smooth oriented geodesically convex spherical curve, and let $-\g$ be its antipodal curve. The curve $\g$ lies in a hemisphere, and $-\g$ lies in the antipodal hemisphere. The spherical belt (topologically, a cylinder) bounded by $\g$ and $-\g$ is the phase space of the outer billiard about $\g$; it is foliated by the arc of the positive tangent great circles to $\g$: these segments have the initial points on $\g$, and the terminal points on $-\g$. This is the vertical foliation that appears in the definition of twist maps.

Likewise, the phase space is foliated by the arc of the negative tangent great circles to $\g$. This makes it possible to define the outer billiard similarly to the planar case: given a point $A$, there is a unique point $x\in\g$ such that the arc of the negative tangent great circle at $x$ contains $A$. The image point $A'$ lies on the  arc of the positive tangent great circle at $x$ at the same spherical distance  from $x$ as point $A$.

The outer billiard map about $-\g$ is conjugated to that about $\g$ by the antipodal involution of the sphere.

\begin{figure}[ht]
\centering
\includegraphics[width=.7\textwidth]{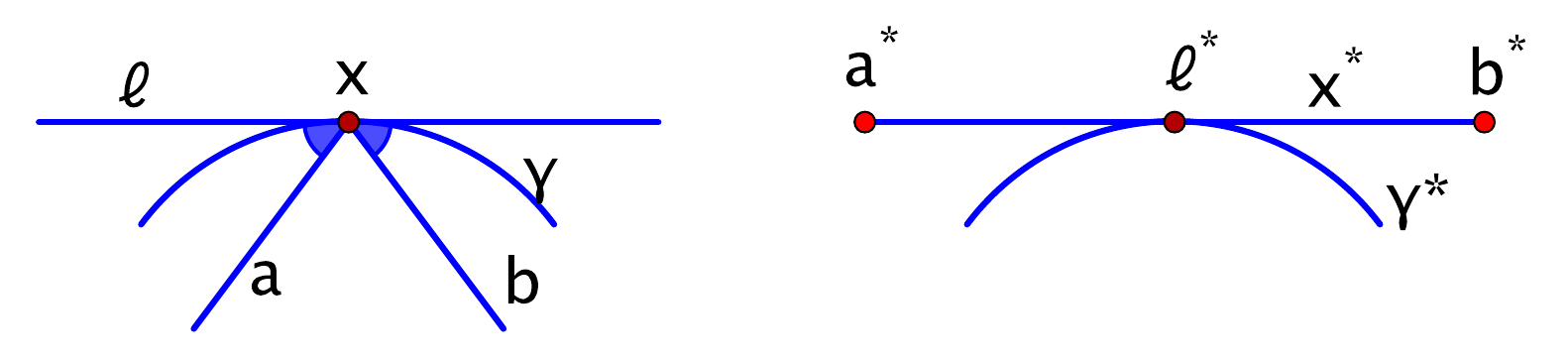}
\caption{The spherical duality conjugates inner and outer billiards.}	
\label{dual}
\end{figure}

Inner and outer billiards in $S^2$ are conjugated by the spherical duality, see Figure \ref{dual}. The spherical duality interchanges oriented great circles with their poles, and the angle between two circles is equal to the spherical distance between their poles. The duality extends to convex smooth curves: the poles of the 1-parameter family of the tangent great circles of a curve $\g$ comprise the dual curve $\g^*$. Equivalently, $\g^*$ is the $\pi/2$-equidistant curve of $\g$, that is, $\g^*$ is the locus of the endpoints of the arcs of the great circles, orthogonal to $\g$ and having length $\pi/2$.

Let $\g$ be a convex smooth curve and $\g^*$ be its dual. Let $L,A,L^*,A^*$ be the perimeters of these curves and the areas of the convex domains bounded by them. Then
$$
L^*=2\pi-A, \ A^*=2\pi-L,
$$
see \cite{An}.

These relations is easy to see if $P$ is a spherical convex $n$-gon and $P^*$ is its dual. Then the angles and the sides lengths of these polygons are related as follows:
$$
\alpha_i^*=\pi-\ell_i,\ \ell_i^*=\pi-\alpha_i.
$$
The Gauss-Bonnet theorem implies
$$
A^*=\sum \alpha_i^* -\pi(n-2)= 2\pi-\sum \ell_i =2\pi-L,
$$
as claimed.  The other equality follows by interchanging $P$ and $P^*$. See, e.g., \cite{FT}, chapter 20, for details.

If $\g$ is a closed smooth convex spherical curve and $P$ is a circumscribed spherical polygon, then the dual polygon $P^*$ is inscribed in the dual curve $\g^*$. It follows that $P$ has the minimal area or perimeter if and only if $P^*$ has the the maximal perimeter or area, respectively.  

See Section 2.1.3 of \cite{AT} for a discussion of symplectic billiards in the spherical and hyperbolic geometries. The convexity of Mather's $\beta$-function implies the spherical version of Dowker's inequalities due to L. Fejes T\'oth \cite{FT2}. 

\subsection{Wire billiards} \label{subsect:wire}

Wire billiards were introduced and studied in \cite{BMT}, see also \cite{Bo1,Bo2,FKM}.

Let $\g(x)$ be a smooth closed arc length parameterized curve in $\R^n$ (a wire). One defines the wire billiard relation in the same way as for the conventional billiards in the plane: chords $\g(x)\g(y)$ and $\g(y)\g(z)$ are in this relation if 
$$
\frac{\partial }{\partial y} (|\g(x)\g(y)|+|\g(y)\g(z)|)=0.
$$
Thus, the phase space for wire billiards is given by oriented chords of $\g$ and the vertical foliation is given by the chords with a fixed initial point.

There is a class of curves, including small $C^2$ perturbations of planar ovals, for which the wire billiard relation is 
a map which then is an area-preserving twist map, called the wire billiard map. This class of curves is given by the following three conditions. 
\begin{enumerate}
\item Any line in $\R^n$ intersects $\g$ in at most two points, and if it intersects at two points, the intersections are with non-zero angles. 
\item The curvature of $\g$ does not vanish.
\item Let $\pi_{xy}$ be the 2-plane spanned by the tangent vector $\g'(x)$ and the chord $\g(y)-\g(x)$. Then for 
every $x,y$ the planes $\pi_{xy}$ and $\pi_{yx}$ are not orthogonal.
\end{enumerate}
If these conditions are satisfied, the generating function of the wire billiard map is given by the same expression as for the conventional planar billiards, $H(x,x')=x'-x-|\g(x)\g(x')|,$
and 
$$
y=-\frac{\partial H}{\partial x}=1-\cos\alpha\in(0,2),
$$
where $\alpha$ is the angle between $\gamma'(x)$ and $\g(x')-\g(x)$. Moreover,
$$
\frac{\partial^2 H(x,x')}{\partial x \partial x'} = -\frac{\cos\varphi\sin \alpha \sin \alpha'}{|\g(x)\g(x')|}<0,
$$
where $\varphi$ is the angle between the planes $\pi_{xy}$ and $\pi_{yx}$.

Let $\g\subset \R^n$ be a curve satisfying the above conditions. Let $R_{p,q}$ be the greatest perimeter of a $q$-gon with the rotation number $p$ inscribed in $\g$, that is, whose vertices lie on $\g$. Then Lemma \ref{lm:ineq} implies that $R_{p,q-1}+R_{p,q+1} < 2R_{p,q}$, a generalization of the Moln\'ar-Eggleston theorem (\ref{eqn:Mol}) for non-planar curves.

\subsection{Minkowski plane} \label{subsect:Mink}\

The Moln\'ar-Eggleston theorem holds in Minkowski planes (2-dimensional Banach spaces), see Theorem 10 in \cite{MSW}. It is geometrically clear that inner and outer length billiards are still twist maps if the Euclidean metric is replaced by a general norm. Therefore, the  Moln\'ar-Eggleston theorem for Minkowski planes should be deducible, as in the Euclidean case, from Lemma \ref{lm:ineq}. We decided not to resolve the details here. 

As for the outer area billiard, we recall that a point and its reflection under the outer area billiard map lie on the same tangent line to the table and having the same distance to the tangent point, see Figure \ref{outer} in Section \ref{subsect:out}. That is, the ratio of the distances is 1. Since the ratio of two distances measured in any norm in $\R^2$ is the same, we see that the reflection rule of outer area billiards is independent of the choice of a norm. This is, of course, no surprise since the generating function is the standard area and does not involve a choice of a norm. 

Finally, the symplectic billiard in the plane simply does not involve a metric in its definition or the reflection rule.

\subsection{Wire symplectic billiards} \label{subsect:wiresymp}

Let $\g(x)$ be a smooth parameterized closed curve in the linear symplectic space $(\R^{2n},\omega)$. We define the symplectic billiard relation on the chords of $\g$ that generalizes symplectic billiards in the plane.

Two chords $\g(x)\g(x')$ and $\g(x')\g(x'')$ of $\g$ are said to be in symplectic billiard relation if $\g(x'')-\g(x) \in T^\omega_{x'} \g$. Here, $T^\omega_{x'} \g:=(T_{x'} \g)^\omega$ is the symplectic orthogonal complement of the tangent line $T_{\g(x')} \g$. Therefore, as for wire billiards, the phase space for symplectic wire billiards is given by oriented chords of $\g$ and the vertical foliation is given by the chords with fixed initial point.

If $x$ and $x''$ are fixed, then $\g(x)\g(x')$ and $\g(x')\g(x'')$ are in symplectic billiard relation if and only if 
$$
\frac{\partial [\omega(\g(x),\g(x'))+\omega(\g(x'),\g(x''))]}{\partial x'} =\omega(\g(x)-\g(x''),\g'(x'))=0,
$$
since $\g'(x')-\g(x) \in T^\omega_{x'} \g$ is equivalent to $\omega(\g(x)-\g(x''),\g'(x'))=0$.

As for wire billiards, this relation does not necessarily define a map. Furthermore, for symplectic wire billiard there is the additional complication that, even if it defines a map, this map need not be a twist map. 

We describe a class of curves $\g$ for which this relation is indeed an area preserving twist map. We recall from \cite{AT1} that a curve $\g$ is called symplectically convex if $\omega(\g'(x),\g''(x))>0$ for all $x$. Consider such a curve, and fix a value $x_0$ of the parameter.  Then the function $F_{x_0} (x) := \omega(\g'(x_0),\g(x))$ has a critical point at $x_0$. Moreover, this zero critical value is a local minimum since
$$
F''_{x_0} (x_0)=\omega(\g'(x_0),\g''(x_0))>0.
$$
The class of curves $\g$ that we consider is given by the following properties:
\begin{enumerate}
\item $\g$ is symplectically convex.
\item For every $x_0$, the function $F_{x_0}: S^1\to \R$ is a perfect Morse function, that is, it has exactly two non-degenerate critical points, a maximum and a minimum.
\end{enumerate}

We call any such curve {\it admissible}. This class of curves is open in the $C^2$-topology, and a sufficiently small perturbation of an oval that lies in a symplectic plane $\R^2 \subset \R^{2n}$ is an admissible curve.

For an admissible curve, we denote by $x^*$ the maximum point of the function $F_{x}$, i.e., $\omega(\g'(x),\g'(x^*))=0$.  We observe that, for an admissible curve $\g$ and a chord $\g(x)\g(x')$, there exists a unique chord $\g(x')\g(x'')$ such that $\g(x'')-\g(x) \in T^\omega_{\g(x')} \g$. 
Indeed, the non-critical level sets of the function $F_{x'}: S^1\to \R$ consist of two points. This makes it possible to extend what we said above about symplectic billiards in the plane to this setting.

Namely, the phase space of wire symplectic billiard is the set of oriented chords $\g(x)\g(x')$ of $\g$ satisfying $x < x' < x^*$ or, equivalently, $\omega(\g(x),\g(x'))>0$. The vertical foliation consists of the chords with a fixed initial point.
As before, we extend the map to the boundary of the phase space by continuity as follows:
$$
(xx) \mapsto (xx),\ (xx^*) \mapsto (x^* x).
$$
The next lemma justifies this statement and repeats a result from \cite{AT}.

\begin{lemma} \label{lm:ok}
Let the wire symplectic billiard map take $xx'$ to $x'x''$. If $x < x' < x^*$, then $x' < x'' < (x')^*$.
\end{lemma}

\proof 
If $x'$ is close to $x$, then $x' < x'' < (x')^*$. 
If $(x')^* < x''$, then, by continuity, we move point $x'$ toward $x$ until $x''=(x')^*$. Then
$$
(xx')\mapsto (x' (x')^*) \mapsto ((x')^* x').
$$
But the map is reversible: if $(xy) \mapsto (yz)$, then $(zy) \mapsto (yx)$. Hence $x=x'$, which is a contradiction.
\proofend

Similarly to the plane case, as the generating function we take
$$
H(x,x')=\int_C \lambda,
$$
where the integral is over the closed curve $C$ made of the arc $\g(x)\g(x')$ of $\g$ and the chord $\g(x') \g(x)$. The curve $C$ is oriented according to the orientation of $\g$ and $\lambda$ is a differential 1-form such that $\omega=d\lambda$. The result does not depend on the choice of such $\lambda$ and is equal to the symplectic area of a surface filling the curve. 

As before, we have 
$$
y=-\frac{\partial H(x,x')}{\partial x}=\frac12\omega(\g'(x),\g(x')-\g(x))
$$
and
$$
\frac{\partial^2 H(x,x')}{\partial x \partial x'} = - \omega(\g'(x),\g'(x')).
$$
Thus, the twist condition is precisely the condition that $\g$ is symplectically convex.

Now consider $q$-gons whose vertices lie on $\g$ and that have the rotation number $p$. The symplectic area of a $q$-gon $(p_1,p_2,\ldots, p_q)$ is 
$$
\frac{1}{2} \sum_{i=1}^q \omega(p_i,p_{i+1}),
$$
where the sum is read cyclically, i.e., $p_{q+1}=p_1$. Let $P_{p,q}$ be the greatest symplectic area of such polygons. Then Lemma \ref{lm:ineq} implies a non-planar generalization of the Dowker theorem:
$P_{p,q-1}+P_{p,q+1} < 2P_{p,q}$.

\subsection{Magnetic billiards} \label{subsect:magn} 

Magnetic billiards were introduced in \cite{RB}; our main reference is \cite{BK}. 

Let $\g$ be a plane oval. 
We assume that the magnetic field has constant strength and is perpendicular to the plane. Then the free path of a charge having a fixed energy is an arc of a circle of radius $R$ (the Larmor radius). When the charge hits the boundary $\g$, it undergoes the billiard reflection, so that the angle of incidence equals the angle of reflection. Unlike the conventional billiards, magnetic billiards is not time-reversible system. It is invariant under simultaneous time and magnetic field reversal, however. Therefore we assume, without loss of generality, that the charge moves in the counterclockwise direction.

Let $k>0$ be the curvature function of $\g$, and $k_{\rm min}$ be its minimal value. The magnetic field is called weak if $1/R < k_{\rm min}$. We consider only the weak field regime here.

\begin{figure}[ht]
\centering
\includegraphics[width=.45\textwidth]{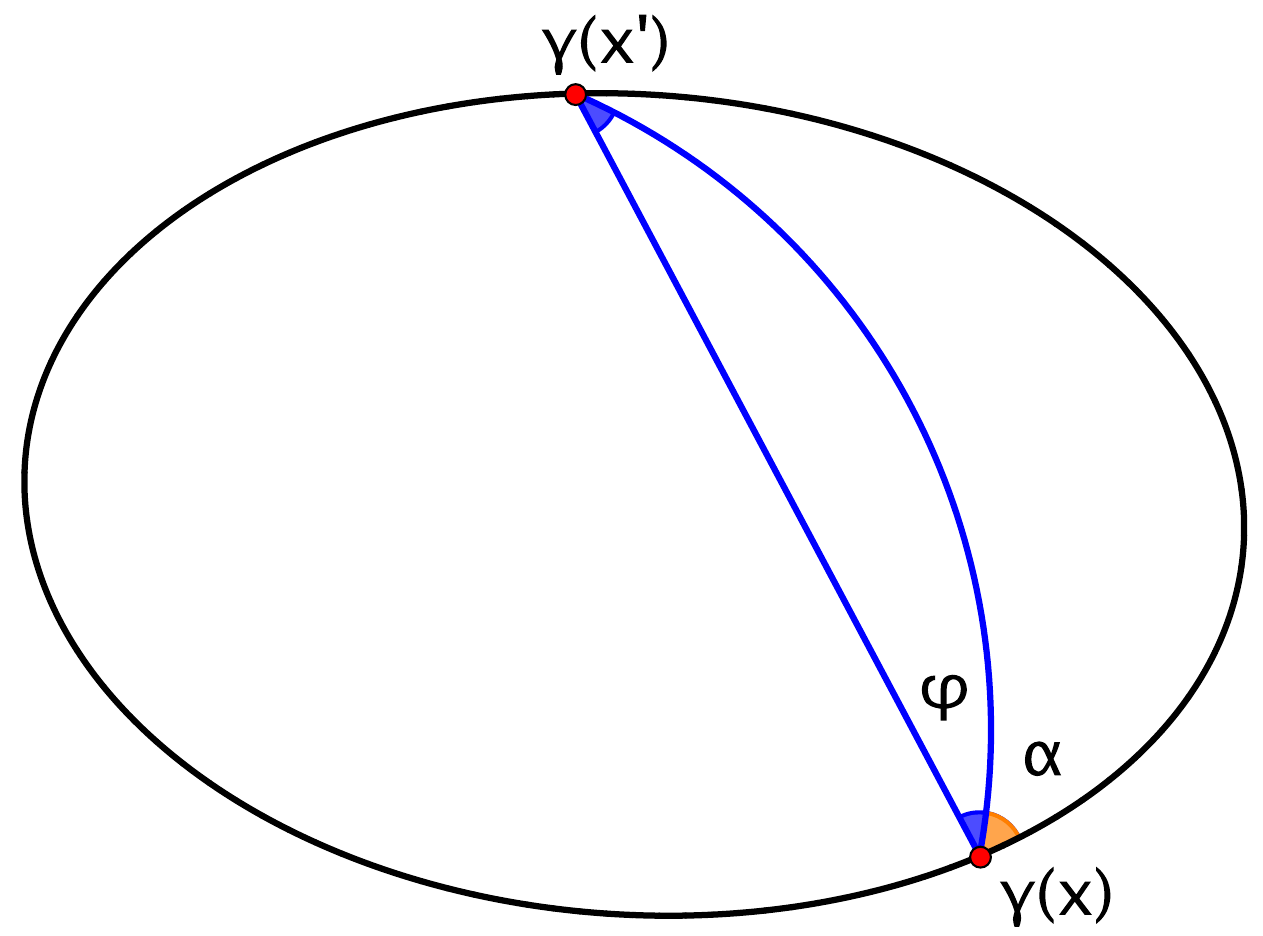}
\caption{Magnetic billiard. Here $\alpha$ is the (orange) angle made by $\g$ and the arc of the circle connecting $\g(x)$ and $\g(x')$. The (blue) angle $\varphi$ is the angle between the arc of the circle and the chord connecting $\g(x)$ and $\g(x')$.}	
\label{magnetic}
\end{figure}

One uses the same coordinates as for conventional billiards: the arc length parameter $x$ on $\g$ and $u=1-\cos\alpha$, where $\alpha$ is the angle made by $\g$ and the trajectory at the starting point. That is, the phase space is the space of oriented chords of $\g$, the vertical foliation consists of the chords with a fixed initial point.

Let $\varphi$ be the angle made by the arc of the trajectory connecting the points $\g(x)$ and $\g(x')$ and the chord connecting these two points, see Figure \ref{magnetic}. 
Then, according to \cite{BK},
$$
\frac{\partial x'}{\partial u} = \frac{|\g(x)\g(x')| \cos\varphi}{\sin\alpha \sin\alpha'}.
$$
Therefore magnetic billiards satisfies the twist condition if the angle $\varphi$ is always acute.

If this assumption holds, then the generating function of the magnetic billiard map is
\begin{equation}\label{eqn:gen_mag_billiard}
H(x,x')=x'-x-\left(\ell + \frac{1}{R} A\right),
\end{equation}
where $\ell$ the the length of the arc $\g(x)\g(x')$ of radius $R$, and $A$ is the area bounded by this arc and the curve $\g$ and lying on the right of this arc  (see again \cite{BK}). When there is no magnetic field, that is, when $R=\infty$, we obtain the generating function of conventional billiard used above. For completeness we also recall from the appendix in \cite{BK}
$$
y=-\frac{\partial H}{\partial x}=1-\cos\alpha\in(0,2),
$$
which is the same as in the conventional billiard case.

We now show that if the magnetic field is weak then the magnetic billiard map is indeed a twist map. The next statement is contained in \cite{BK} as Lemma D1 and the following Corollary. We provide a slightly different proof here for convenience.

\begin{lemma} \label{lm:weak}
If the magnetic field is weak, magnetic billiard map is a twist map. More precisely, if $1/R < k_{\rm min}$, then $\varphi < \pi/2$.
\end{lemma}

\begin{figure}[ht]
\centering
\includegraphics[width=.3\textwidth]{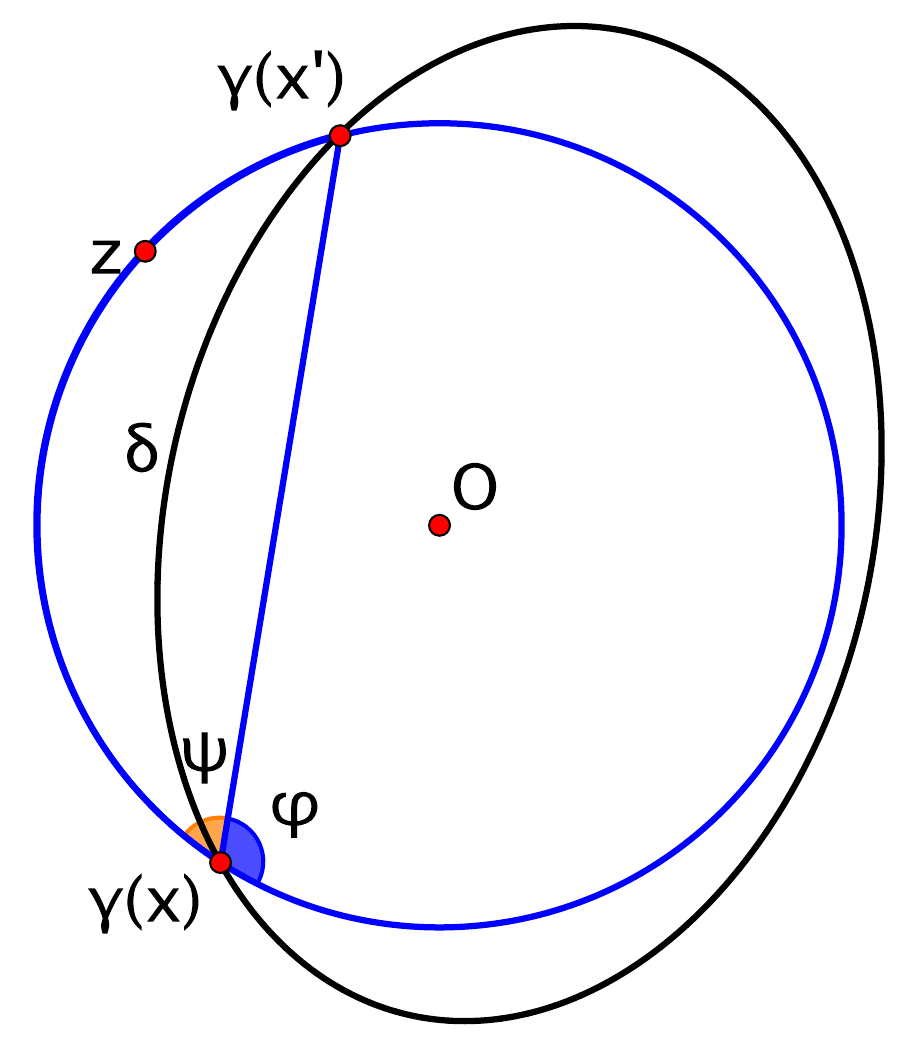}
\caption{Concerning Lemma \ref{lm:weak}.}	
\label{Schur}
\end{figure}

\proof
Let $\psi=\pi-\phi$ be the complementary angle, see Figure \ref{Schur}. We shall prove the equivalent statement: if $\psi$ is acute, then there exists a point on the arc $\delta$ of the curve $\g$ from $\g(x)$ to $\g(x')$ where the curvature $k \le 1/R$.

We use two  geometric facts. The first is a lemma due to Schur, see, e.g., \cite{Gu}, asserting the following. Let $\g_1(s)$ and $\g_2(s)$ be two smooth convex arc length parameterized curves of the same length, such that $k_1(s) < k_2(s)$ for all $s$. Then the chord subtended by $\g_1$ is greater than that the chord subtended by $\g_2$ for the same interval of the parameter.

The second is the following lemma. If a closed convex curve $\g_2$ lies inside (the domain bounded by) the closed convex curve $\g_1$, then the length of $\g_1$ is greater than or equal to the length of $\g_2$. This is an easy consequence of the Crofton formula, see, e.g., \cite{FT}, Lecture 19, and it does not require the curves to be smooth, nor for $\g_2$ to be contained strictly inside of $\g_1$.

Now, consider the clockwise oriented arcs of the curve $\g$ and of the circle starting at the point $\g(x)$. Consider the first intersection point of these arcs. In Figure \ref{Schur}, this point is $\g(x')$ but, in general, it may lie closer to point $\g(x)$ on both curves. In order  not to complicate the notation, we assume that indeed the first intersection point is $\g(x')$.

Then we have two nested convex closed curves: the first is made by the arc of the circle and the chord $\g(x')\g(x)$, and the second by the arc $\delta$ and the same chord. By assumption the second curve lies inside the first. It follows that the arc of the circle has length greater than or equal to that of $\delta$. 

Now we argue by contradiction and assume $1/R < k$. Let $z$ be the point on the arc of the circle such that $\g(x) z$ has the same length as $\delta$. Since $1/R < k$, Schur's Lemma implies that
$|\g(x) z| > |\g(x')\g(x)|$. But the assumption that $\psi < \pi/2$ implies that the center of the circle lies on the right of the chord $\g(x')\g(x)$, therefore $|\g(x')z|<|\g(x')\g(x)|$, cf.~Figure \ref{Schur}. This is a contradiction.
\proofend

\begin{remark}
{\rm The same argument proves that if $1/R < k_{\rm min}$, then a circle of radius $R$ intersects the curve $\g$ in at most two points.
}
\end{remark}

As before, we apply Lemma \ref{lm:ineq} to obtain the following result. Consider curvilinear $q$-gons with rotation number $p$ inscribed in the oval $\g$ whose oriented sides are counterclockwise arcs of radius $R$ with $1/R < k_{\rm min}$. Let $M_{p,q}$ be the greatest value of 
$$
P-\frac1R A
$$
where $P$ is the perimeter of the curvilinear $q$-gon and $A$ in the area \emph{enclosed} by the curvilinear polygon. Then one has the Dowker-style inequality: $M_{p,q-1}+M_{p,q+1} < 2M_{p,q}$.

Note that the area term in the generating function \eqref{eqn:gen_mag_billiard} is the area between the curve $\g$ and the respective arc of the polygon. Since the polygon is closed the area between the polygon and the curve and the area enclosed by the polygon differ by the total area enclosed by $\g$. This does not change the inequalities.

\subsection{Outer magnetic billiards} \label{subsect:outmagn} 

\begin{figure}[ht]
\centering
\includegraphics[width=.4\textwidth]{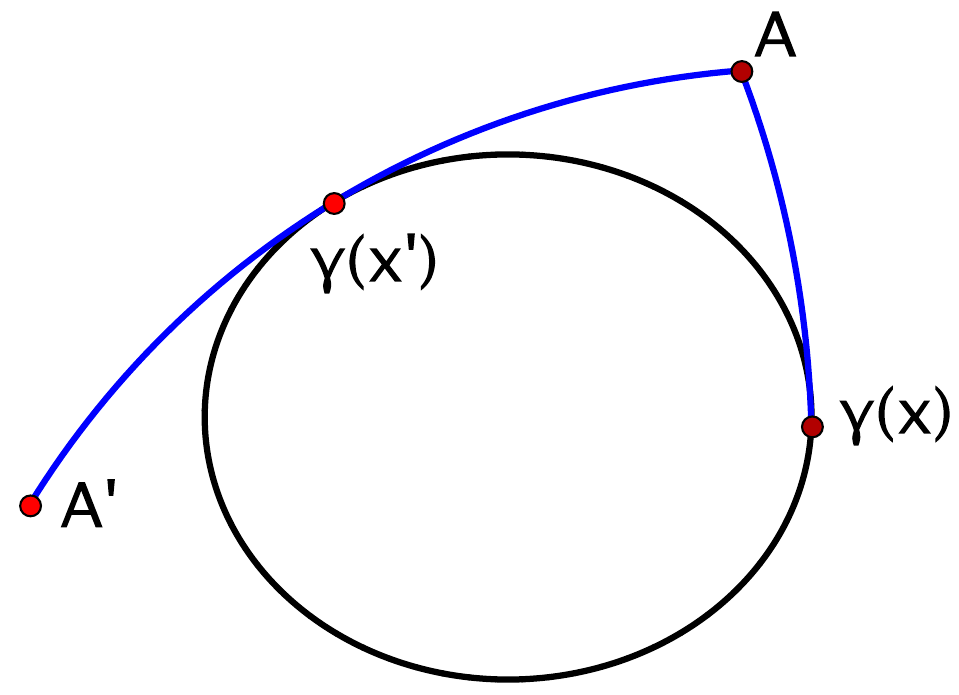}\qquad\quad
\includegraphics[width=.35\textwidth]{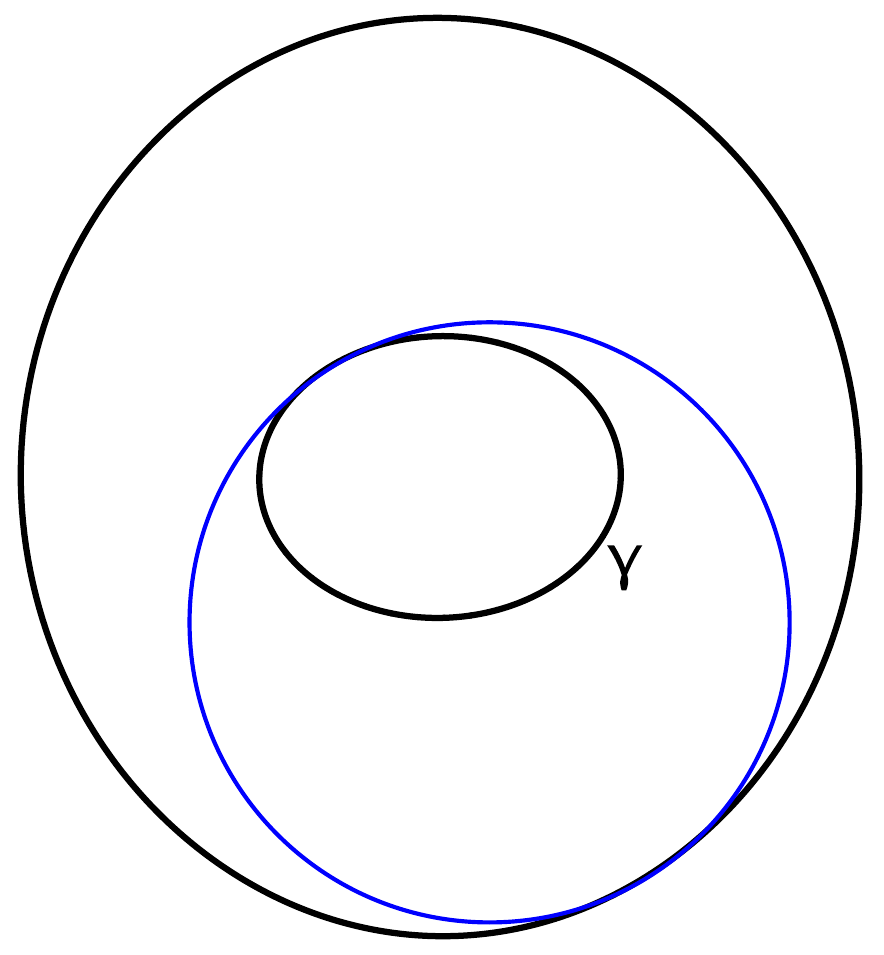}
\caption{Left: outer magnetic billiard map; right: its annulus of definition.}	
\label{outmag}
\end{figure}

The outer magnetic billiard map is defined similarly to the outer billiard map, see Section \ref{subsect:out}, but the tangent lines are replaced by the tangent arcs of circles of a fixed (Larmor) radius, greater than the greatest radius of curvature of the ``billiard" curve $\gamma$, see Figure \ref{outmag} on the left. That is, we are again in the weak magnetic field regime. We assume that $\gamma$ and the circles are positively oriented and that the orientations agree at the tangency points. 

Unlike the outer billiard map, this map is only defined in an annulus whose inner boundary is $\gamma$ and the outer boundary is the envelope of the Larmor circles tangent to $\gamma$ which, by Huygens' principle, is equidistant from $\gamma$,  see Figure \ref{outmag} on the right.

This map is area preserving, and the value of the generating function at the point $A$ is the area of the curvilinear triangle bounded by the two Larmor arcs through $A$ and the curve $\gamma$. The vertical foliation is given by the forward Larmor half-circles tangent to $\gamma$.

Outer magnetic billiard was introduced in \cite{BM}, where it is shown that this system is isomorphic to magnetic billiard. The correspondence between the latter and the former is given by the map that assigns to the arcs of Larmor circles (inside the billiard table) their centers. The resulting outer magnetic billiard curve is equidistant to the magnetic billiard curve.  

A Dowker-style geometric inequality results. Let $\gamma$ be an oval, and let 
$N_{p,q}$ denote the minimal area of a circumscribed curvilinear $q$-gon with rotation number $p$, whose sides are arcs of a fixed radius greater that the greatest curvature radius of $\gamma$. Then $N_{p,q-1}+N_{p,q+1} > 2N_{p,q}$. We note that, unlike for inner magnetic billiards, in this case the geometric inequality involves areas only and not a combination of areas and lengths.  

\subsection{Remarks on the theory of interpolating Hamiltonians}

We conclude with remarks and questions concerning an application of the theory of interpolating Hamiltonians \cite{MM,Mel,PS} in 
convex geometry, namely, to approximation of smooth convex curves by polygons.

It is proved in \cite{MM} that the maximal perimeter of the inscribed $q$-gons 
$R_{1,q}$ has an asymptotic expansion as $q\to \infty$:
$$
R_{1,q} \sim {\rm Perimeter} (\g) + \sum_{k=1}^\infty \frac{c_k}{q^{2k}},
$$
with the first non-trivial coefficient 
$$
c_1=-\frac{1}{24} \left(\int_{\g} k^{\frac23} ds\right)^3,
$$
where $k(s)$ is the curvature of $\g$ and $s$ is the arc length parameter. If we read Dowker's inequality as
$$
R_{1,q+1}-R_{1,q}\leq R_{1,q}-R_{1,q-1},
$$
we see that the successive approximation by polygons always improves. 

Also the limiting $q\to\infty$ distribution of the vertices of the approximating inscribed $q$-gons is uniform with respect to the density  $k^{\frac23} ds$.
Likewise for other rotation numbers.

Similarly, one has for the minimal area of the circumscribed $q$-gons
$$
Q_{1,q} \sim {\rm Area}(\g) + \sum_{k=1}^\infty \frac{c_k}{q^{2k}},
$$
with 
$$
c_1=\frac{1}{24} \left(\int_{\g} k^{\frac13} ds\right)^3,
$$
see \cite{Ta95}. For the maximal area of the inscribed $q$-gons, one has
$$
P_{1,q} \sim {\rm Area}(\g) + \sum_{k=1}^\infty \frac{c_k}{q^{2k}},
$$
with 
$$
c_1=-\frac{1}{12} \left(\int_{\g} k^{\frac13} ds\right)^3,
$$
see \cite{AT} and \cite{MV}.

In the last two cases, the limiting  $q\to\infty$ distribution of the vertices of the approximating $q$-gons on the curve $\g$ is uniform with respect to the density  $k^{\frac13} ds$, that is, uniform with respect to the affine length parameter.

It would be interesting to find explicit expressions for the coefficients $c_k$ in these formulas  (see \cite{Lu,MM} for $c_1$ and $c_2$). A closely related problem of describing the coefficients of the expansion of Mather's $\beta$-function at zero for Birkhoff billiards is addressed in \cite{So} and, for ellipses, in \cite{Bi}, and for symplectic and outer billiards, very recently, in \cite{BBN}.

There are  other ways to measure the quality of approximation of a convex curve by polygons, for example,  one can use the area of the symmetric difference of an oval and an approximating polygon, or its analog, replacing area by perimeter. A wealth of results in this directions is available, see \cite{Bo,Fo,Gr}. Are these results related to the $\beta$-function  and the interpolating Hamiltonians of some billiard-like dynamical systems?

\bigskip

{\bf Acknowledgements}. 
PA acknowledge funding by the Deutsche Forschungsgemeinschaft (DFG, German Research Foundation) through Germany's Excellence Strategy EXC-2181/1 - 390900948 (the Heidelberg STRUCTURES Excellence Cluster), the Transregional Colloborative Research Center CRC/TRR 191 (281071066). ST was supported by NSF grant DMS-2005444 and by a Mercator fellowship within the CRC/TRR 191, and he thanks Heidelberg University for its invariable hospitality.

\end{document}